\documentclass{amsart}

\begin{document}

\title{Arithmetic analogues of Hamiltonian systems}
\bigskip

\def \h{\hat{\ }}
\def \cO{\mathcal O}
\def \ra{\rightarrow}
\def \bZ{{\mathbb Z}}
\def \cP{\mathcal V}
\def \cH{{\mathcal H}}
\def \cB{{\mathcal B}}
\def \d{\delta}
\def \cC{{\mathcal C}}
\def \jor{\text{jor}}

\newtheorem{THM}{{\!}}[section]
\newtheorem{THMX}{{\!}}
\renewcommand{\theTHMX}{}
\newtheorem{theorem}{Theorem}[section]
\newtheorem{corollary}[theorem]{Corollary}
\newtheorem{lemma}[theorem]{Lemma}
\newtheorem{proposition}[theorem]{Proposition}
\newtheorem{thm}[theorem]{Theorem}
\theoremstyle{definition}
\newtheorem{definition}[theorem]{Definition}
\theoremstyle{remark}
\newtheorem{remark}[theorem]{Remark}
\newtheorem{example}[theorem]{\bf Example}
\numberwithin{equation}{section}
\address{Department of Mathematics and Statistics\\University of New Mexico \\
  Albuquerque, NM 87131, USA\\}
\email{buium@math.unm.edu} 
\maketitle

\bigskip

\medskip
\centerline{\bf Alexandru Buium}
\bigskip

\centerline{\it Dedicated to Emma Previato's $65$th birthday}

\bigskip

\begin{abstract} 
The paper reviews various arithmetic analogues of Hamiltonian systems 
introduced in \cite{foundations, gln2, BYM, euler, canonical}, and presents some new facts suggesting ways to relate/unify these examples.
\end{abstract}

\section{Introduction}

\subsection{Aim of the paper}
In a series of papers starting with \cite{char} an arithmetic analogue of the concept of (not necessarily linear) ordinary differential equation was developed; cf. the monographs \cite{book} and \cite{foundations}
for an exposition of (and references for) this theory and \cite{pj, dmf, local, fourier} for some purely arithmetic applications. 
(There is a version for partial differential equations \cite{laplace, pde1, pde2, foundations} which we will not discuss here.)
The rough idea behind this theory is to replace differentiation operators 
$$y\mapsto \d_x y:=\frac{dy}{dx},$$
acting on smooth functions $y=y(x)$ in one variable $x$, by {\it Fermat quotient operators} with respect to an odd rational prime $p$,
$$z\mapsto \d_p z:=\frac{\phi_p(z)-z^p}{p},$$ 
acting on  
a complete discrete valuation ring $A$ with maximal ideal generated by $p$ and perfect residue field; 
here we denoted by $\phi_p$  the unique Frobenius lift on $A$, i.e. the unique ring endomorphism of $A$ whose reduction mod $p$ is the $p$-power Frobenius. 
Then classical  differential equations,
$$F(y, \d_x y,...,\d_x^n y)=0,$$
where $F$ is a smooth function,
are replaced by {\it arithmetic differential equations},
$$F(z,\d_p z,...,\d_p^n z)=0,$$
where $F$ is a $p$-adic limit of polynomials with coefficients in $A$. More generally one can consider systems of such equations.
One would then like to introduce arithmetic analogues of some remarkable classical ordinary differential equations such as: 

\medskip

\noindent {\it  Linear differential equations, 
 Riccati equations, 
 Weierstrass equation satisfied by elliptic functions, 
Painlev\'{e} equations, 
Schwarzian equations satisfied by modular forms, 
Ramanujan differential equations satisfied by  Einsenstein series, 
Euler equations for the rigid body, 
 Lax equations, etc.}

\medskip

\noindent A first temptation, in developing the theory, would be to consider  the function(s) $F$ in each of these classical equations and
formally replace, in the corresponding equation, the quantities $\d_x^i y$ by the quantities $\d_p^i z$. This strategy of preserving $F$, i.e., the {\it shape of the equations}, turns out to destroy the underlying {\it geometry of the equations} and, therefore, seems to lead to a dead end.  Instead, what turns out to lead to an interesting theory (and to interesting applications) is to discover how to change the shape of the equations (i.e. change $F$) so that the geometry of the equations is (in some sense) preserved. The first step in then to ``geometrize" the situation by introducing {\it arithmetic jet spaces} \cite{char} which are arithmetic analogues of the classical jet spaces of Lie and Cartan and to seek a conceptual approach towards arithmetic analogues of the classical equations we just mentioned; this has been done in the series of papers and in the monographs
mentioned in the beginning of this paper, for the list of classical equations we just mentioned. In achieving this one encounters various obstacles and surprises. For instance, the question, 
{\it ``What is the  arithmetic analogue of  linear differential equations?"}  
 is already  rather subtle; indeed, linearity of arithmetic differential equations turns out to be not an absolute but, rather, a relative concept; more precisely there is no concept of linearity for {\it one} 
arithmetic differential equation but there is a concept of an arithmetic differential equation being {\it linear} with respect to another arithmetic differential equation. We refer to \cite{foundations, gln2, gln3} for the problem of linearity. 
The question we would like to address in this paper is a different one (although a related one), namely:

\medskip

{\it ``What is the  arithmetic analogue of a Hamiltonian system?"}  

\medskip

\noindent A number of remarkable classical differential equations admit  Hamiltonian structures; this is
the case, for instance, with the following 3 examples: Painlev\'{e} VI equations, Euler equations, and (certain special) Lax equations.
Arithmetic analogues of these 3 types of equations have been developed, in 3 separate frameworks, in a series of papers as follows: the Painlev\'{e} case in \cite{BYM};
the Euler case in \cite{euler, canonical}; and the Lax case in \cite{gln2}, respectively. Cf. also \cite{foundations}. One is tempted to believe that these 3 examples are pieces of a larger puzzle.
The aim of the present paper is to review these 3 examples and attempt to give hints as to a possible unification of these 3 pictures by proving  some new facts (and providing some new  comments) that connect some of   the dots. The task of setting up a general framework (and a more general array of examples) for an arithmetic Hamiltonian formalism is still elusive; we hope that the present paper contains clues as to what this general formalism could be.

\subsection{Structure of the paper}
Each of the following sections contains $2$ subsections. In each section the first of the subsections
offers  a treatment of the classical differential setting while the second subsection offers a treatment of the arithmetic differential setting. 
Section 2 of the paper is devoted to an exposition of the main general concepts. Sections 3, 4, 5 are devoted to the main examples under considerations: the Painlev\'{e} equations, the Euler equations, and the Lax equations respectively. 

The main new results of the paper are Theorems \ref{new1} and \ref{new2}
 which hint
towards a  link between the formalisms in the Painlev\'{e} and  Euler examples. 

\subsection{Ackowledegment} 
The present work was partially supported  by the Institut des Hautes \'{E}tudes Scientifiques in Bures sur Yvette, 
and by  the Simons Foundation (award 311773). The author would also like to express his debt
 to Emma Previato for an inspiring collaboration and a  continuing interaction.

\section{General concepts}

\subsection{The classical case}
We begin by reviewing some of the main concepts in the theory of classical (ordinary) differential equations.
We are interested in the purely algebraic aspects of the theory so we will place ourselves in the context  of differential algebra \cite{Kolchin}; our exposition  follows \cite{ajm93, foundations, BYM}. 

Let us start with a  ring $A$ equipped with a derivation $\d^A:A\ra A$ (i.e., an additive map satisfying the Leibniz  rule).  For simplicity we assume $A$ is Noetherian. Let $X$ be a scheme of finite type over $A$. One defines the {\it jet spaces} $J^n(X)$ of $X$ (over $A$) as follows. If $X$ is affine,
\begin{equation}
\label{X}
X=Spec\ A[x]/(f),\end{equation}
with $x$ an $N$-tuple of indeterminates and $f$ a tuple of polynomials, then one sets
\begin{equation}
\label{J}
J^n(X):=Spec\ A[x,x',...,x^{(n)}]/(f,\d f, ... ,\d^nf)\end{equation}
where $x',...,x^{(n)}$ are new $N$-tuples of  indeterminates and $\d=\d^{\text{univ}}$ is the unique (``universal") derivation on the ring of polynomials in infinitely many indeterminates,
$$A[x,x',...,x^{(n)},...],$$
 extending $\d^A$, and satisfying 
\begin{equation}
\d x=x',...,\d x^{(n)}=x^{(n+1)},... \end{equation}
One gets induced derivations
$$\d=\d^{\text{univ}}:\cO(J^n(X))\ra \cO(J^{n+1}(X)).$$
If $X$ is not necessarily affine then one defines $J^n(X)$ by gluing $J^n(X_i)$ where $X=\cup X_i$ is a Zariski affine open cover. The family $(J^n(X))_{n\geq 0}$ has the structure of  a projective system of schemes depending functorially on $X$, with $J^0(X)=X$. 
If $X$ is smooth and descends to (a scheme over)  
{\it the ring of $\d$-constants} of $A$,
$$A^{\d}=\{a\in A;\ \d a=0\},$$
  then $J^1(X)$ identifies with the (total space of the) tangent bundle $T(X)$ of $X$; if we drop the condition that $X$ descend to $A^{\d}$ then $J^1(X)$ is only a torsor under $T(X)$.

If $G$ is a group scheme over $A$ then $(J^n(G))_{n\geq 0}$ forms a projective system of group schemes over $A$; if $A$ is a field of characteristic zero, say,  the kernels of the homomorphisms $J^n(G)\ra J^{n-1}(G)$
are isomorphic as algebraic groups to powers of the additive group scheme ${\mathbb G}_a$.

By a  {\it differential equation} on $X$ we understand a closed subscheme of some $J^n(X)$.
By a $\d^A$-{\it flow}
on $X$ we will understand a derivation $\d^X$ on the structure sheaf of $X$, extending $\d^A$; giving a $\d^A$-flow is equivalent to giving a section of the canonical projection $J^1(X)\ra X$, and hence, to giving a differential equation $Z\subset J^1(X)$ for which the projection $Z\ra X$ is an isomorphism. A {\it prime integral} for a $\d^A$-flow $\d^X$ is a function $H\in \cO(X)$ such that $\d^XH=0$.
For any $A$-point 
$$P\in X(A),\ \ \ P:Spec\ A \ra X,$$
 one defines the jets 
$$J^n(P)\in J^n(X)(A),\ \ \ J^n(P):Spec\ A\ra J^n(X),$$
 as the unique morphisms lifting $P$ that are compatible with the actions of $\d^{\text{univ}}$
 and $\d^A$. 
A {\it solution} (in $A$) for a differential equation $Z\subset J^n(X)$ is an $A$-point 
$P\in X(A)$ such that $J^n(P)$ factors through $Z$. If $P$ is a solution to the differential equation defined by a $\d^A$-flow and if $H$ is a prime integral for that $\d^A$-flow
then $\d^A(H(P))=0$; intuitively $H$ is ``constant" along any solution. If $X$ is affine 
and $Z\subset J^1(X)$
is the differential equation corresponding to a $\d^A$-flow $\d^X$ on $X$ then a point $P\in X(A)$ is a solution to $Z$ if and only if the ring homomorphism $P^*:\cO(X)\ra A$ defined by $P$ 
satisfies $\d^A\circ P^*=P^*\circ \d^X$.

For any smooth affine scheme $Spec\ B$ over $A$ we may consider  the algebraic deRham complex of abelian groups,
\begin{equation}
B\stackrel{d}{\longrightarrow} \Omega^1_{B/A} \stackrel{d}{\longrightarrow} \Omega^2_{B/A}\stackrel{d}{\longrightarrow} ...,
\end{equation}
where $\Omega^i_{B/A}:=\wedge^i \Omega_{B/A}$,  $\Omega_{B/A}$ the (locally free) $B$-module of K\"{a}hler differentials.

Recall that for any $A$-algebra $B$ a {\it Poisson structure/bracket} on $B$ (or on $Spec\ B$) is an $A$-bilinear map
$$\{\ ,\ \}:B\times B\ra B$$
which is a derivation in each of the arguments and defines a structure of Lie $A$-algebra on $B$.

In what follows we would like to review the classical concept of {\it Hamiltonian} system/equation from a purely algebro-geometric viewpoint; later we will introduce its arithmetic analogues. We will restrict ourselves to the case of affine surfaces since our main examples fit into this setting.

Let $S=Spec \ B$ be a smooth affine surface (i.e. smooth   scheme of relative dimension $2$) over a Noetherian ring $A$. 

By a {\it symplectic form} we will understand a basis $\eta$ of the $B$-module $\Omega^2_{B/A}$. (The usual condition that this form be closed is automatically satisfied because $S$ is a surface.)
By a {\it contact form} we will understand an element  $\nu\in \Omega^1_{B/A}$ such that
$d\nu$ is a symplectic form.
Given a symplectic form $\eta$ on $S=Spec\ B$ one can define a Poisson structure on $S$ by the formula:
$$\{ \ ,\ \}_{\eta}:B\times B\ra B,\ \ \ \{ f, g\}_{\eta}:=\frac{df \wedge dg}{\eta}.$$

Assume now 
that we are given a derivation $\d^A:A\ra A$. 
 Recall that by a {\it $\d^A$-flow} on $S$ we mean a derivation $\d:=\d^B:B\ra B$ 
extending our derivation $\d^A:A\ra A$. Recall that $\d^B$ induces then unique additive maps, referred to as {\it Lie derivatives},
$$\d=\d^{\Omega^i}:\Omega^i_{B/A}\ra \Omega^i_{B/A},\ \ \ i=0,1,2,$$
such that $\d^{\Omega^0}=\d^B$, 
$\d$ commute with $d$,  and $\d$  induce a derivation on the exterior algebra $\wedge \Omega_{B/A}$. 

Recall  that a function $H\in B$ is called a {\it prime integral} for $\d^B$ (or a $\d^B$-{\it constant}) if $\d^BH=0$.

Let us say that a $\d^A$-flow $\d=\d^B$ on $S$ is {\it Hamiltonian} 
(or, more accurately,  {\it symplectic-Hamiltonian})
with respect to a symplectic form $\eta$ on $S$ if
\begin{equation}
\label{Lieder}
\d\eta=0.\end{equation}

Let us say that 
$\d^A$-flow $\d=\d^B$ on $S$ is {\it Hamiltonian} (or, more accurately,  {\it Poisson-Hamiltonian}) 
with respect to a Poisson structure $\{\ ,\ \}$ on $B=\cO(S)$ 
if $S$ descends to a smooth scheme $S_0=Spec\ B_0$ over $A_0:=A^{\d}$, with 
$\{B_0,B_0\}\subset B_0$,
 and 
 there exists a function (called {\it Hamiltonian}) $H\in B_0$ such that
\begin{equation}
\label{poisson}
\d f=\{f,H\},\ \ \text{for all}\ \ f\in B_0.\end{equation}
A direct computation with \'{e}tale coordinates shows that if a $\d^A$-flow $\d=\d^B$ on $S$ is {\it Hamiltonian} with respect to the Poisson structure 
$\{\ ,\ \}_{\eta}$ on $B$ attached to a symplectic form $\eta$ on $S$ that comes from $S_0$ then
$\d$ is Hamiltonian with respect to $\eta$; moreover if  $H$ is a Hamiltonian then, trivially,  $H$ is a prime integral for $\d$. 
As we shall see the examples of Painlev\'{e} and Euler equations are symplectic-Hamiltonian   but {\it not} Poisson-Hamiltonian simply because the surfaces $S$ on which these equations ``live" do not descend to surfaces $S_0$ over the constants. A large class of  examples coming from Lax equations are Poisson-Hamiltonian.
Both symplectic-Hamiltonian and Poisson-Hamiltonian equations have arithmetic analogues.

The discussion above has, of course, a higher dimensional analogue in which $S$ is a smooth affine scheme of arbitrary dimension; the Poisson-Hamiltonian picture is valid word for word; the symplectic-Hamiltonian picture has to be modified  by asking that $S$ have relative dimension $2d$, $\eta$ be a closed $2$-form, and $\eta^d$ be invertible.

Going back to the case when $S$ is a surface, let $\nu$ be a contact form, let $\d=\d^S$ be a $\d^A$-flow on $S$, consider the symplectic form
$\eta:=d\nu$, and define the {\it Euler-Lagrange} form
$$\epsilon:=\d\nu\in \Omega^1_{B/A}.$$
Since $d$ and $\d$ commute we have that $\d^S$ is Hamiltonian with respect to $\eta$ if and only if $\epsilon$ is closed, i.e. $d\epsilon=0$.
If in addition $\epsilon$ is exact, i.e. $\epsilon=d{\mathcal L}$ for some ${\mathcal L}\in \cO(S)$, we call ${\mathcal L}$ a {\it Lagrangian} for $(\nu,\d^S)$.

A special case  that plays a role in the theory is that in which our surface $S$ is the first jet space of a smooth curve $Y$ over $A$,
$$S=J^1(Y).$$ 
In this case a $1$-form $\nu$ on $S$ is called {\it canonical} if $\nu=f\beta$
where $f\in B=\cO(S)$ and $\beta$ is a pull-back of a $1$-form on $Y$.
Assume $\nu$ is a canonical contact form, assume $\d=\d^S$ is Hamiltonian with respect to $\eta:=d\nu$,  assume  $\epsilon$ is exact with Lagrangian ${\mathcal L}$, and assume $x\in \cO(Y)$ is an \'{e}tale coordinate on $Y$. Assume in addition that $x,\d x$ are \'{e}tale coordinates on $S$ (which is ``generically the case" and is automatic, for instance,  for the {\it canonical} $\d^A$-flows to be introduced below). Then there are unique $A$-derivations
$\frac{\partial}{\partial x}, \frac{\partial}{\partial \d x}$ on $\cO(S)$
sending $x,\d x$ into $1,0$ and $0,1$ respectively. It is then trivial to check that
$\d \left( \frac{\partial {\mathcal L}}{\partial \d x}\right)=\frac{\partial {\mathcal L}}{\partial x}$ in $\cO(S)$.
 In particular if $Z\subset J^1(S)$ is the differential equation corresponding to the $\d^A$-flow $\d^S$ on $S$ then any solution in $S(A)$ to $Z$ 
 is a solution to the {\it Euler-Lagrange  equation} $EL(Z) \subset J^1(S)$ defined  by
$$\d^{\text{univ}} \left( \frac{\partial {\mathcal L}}{\partial \d x}\right)-\frac{\partial {\mathcal L}}{\partial x}\in \cO(J^1(S)).$$

Contact forms $\nu$ that are canonical should be viewed as generalizing the {\it canonical contact forms} on  cotangent bundles in differential geometry (and classical mechanics); also our Lagrangians and Euler-Lagrange equation correspond, formally, to the Lagrangians and Euler-Lagrange equation
in classical mechanics.  Note however the following discrepancy with the usual definition in differential geometry: our $J^1(Y)$ is related to (is a torsor 
under) the {\it tangent} bundle while in classical differential geometry {\it canonical forms} live on the  {\it cotangent} bundle. This discrepancy is resolved, in usual differential geometry, by identifying the tangent and the cotangent bundle via $d\nu$; in our setting (when $J^1(Y)$ is not a trivial torsor) no such identification is available.
By the way, as we shall explain,  it is the definition of {\it canonical contact form}  that we just gave above 
(and not the usual definition  in differential geometry)  that will have an arithmetic analogue.

Finally, we make the following definition: 
 a {\it canonical $\d^A$-flow} on $S=J^1(Y)$ is  a $\d^A$-flow $\d^{J^1(Y)}$ on 
 $J^1(Y)$
  with the property that the composition of 
  $$\d^{J^1(Y)}:\cO(J^1(Y))\to \cO(J^1(Y))$$ 
with the pull back map 
$$\cO(Y)\to \cO(J^1(Y))$$ equals the universal
derivation 
$$\d^{\text{univ}}:\cO(Y)\to \cO(J^1(Y)).$$
 By the way, notice  that one has a natural
closed embedding 
$$\iota:J^2(Y)\to J^1(J^1(Y)).$$
Then one checks that a $\d^A$-flow $\d^{J^1(Y)}$ on $J^1(Y)$ is canonical if and only if
the section $J^1(Y)\ra J^1(J^1(Y))$ defined by $\d^{J^1(Y)}$ factors through $\iota$. Also notice that if $x\in \cO(Y)$ is an \'{e}tale coordinate and $\d=\d^{J^1(Y)}$ is a canonical flow then $x,\d x$ are \'{e}tale coordinates on $S$.
The  concept of canonical $\d^A$-flow is an algebraic version of a 
classical concept related to second order ODEs (for instance Painlev\'{e} equations)
and has an arithmetic analogue.

\subsection{The arithmetic case}

Let  $p$ be a rational odd prime.  If $B$ is a ring a {\it Frobenius lift} on $B$  is a ring endomorphism $\phi=\phi^B:B\ra B$ whose reduction mod $p$ is the $p$-power Frobenius on $B/pB$. Similarly if $X$ is a scheme or a $p$-adic formal scheme a {\it Frobenius lift} on $X$  is an endomorphism $\phi=\phi^X:X\ra X$ whose reduction mod $p$ is the $p$-power Frobenius on the reduction of $X$ mod $p$. 
Let $A$ be a complete discrete valuation ring with maximal ideal generated by $p$ and perfect residue field $k=A/pA$; we fix this $A$ once and for all in the discussion below. Such an $A$ is uniquely determined up to isomorphism by $k$ and possesses a unique Frobenius lift $\phi=\phi^A:A\ra A$.
For any $A$-algebra $B$ and any scheme or $p$-adic formal scheme $X$ over $A$ Frobenius lifts on $B$ or $X$ will be tacitly assumed  to be compatible with the Frobenius lift on $A$. For any Noetherian $A$-algebra $B$ and Noetherian scheme $X$ over $A$ we denote by $\widehat{B}$ and $\widehat{X}$ 
the $p$-adic completions of $B$ and $X$ respectively. We also define the K\"{a}hler differentials on the formal scheme $\widehat{X}$ by
\begin{equation}
\label{maxim}
\Omega_{\widehat{X}}=\lim_{\leftarrow} \Omega_{X_n/A_n}\end{equation}
where $A_n=A/p^nA$, $X_n=X\otimes A_n$. If $X$ is smooth over $A$ and $\phi$ is a Frobenius lift on $\widehat{X}$ then $\phi$ naturally induces  additive maps
\begin{equation}
\label{phipep}
\frac{\phi^*}{p^i}:\Omega^i_{\widehat{X}}\ra \Omega^i_{\widehat{X}},\end{equation}
where 
$\Omega^i_{\widehat{X}}:=\wedge^i \Omega_{\widehat{X}}$. 

Given a ring $B$ which is $p$-torsion free (i.e., $p$ is a non-zero divisor in $B$) a map of sets 
$$\d=\d^B:B\ra B$$ will be called a $p$-{\it derivation} if the map
$$\phi=\phi^B:B\ra B, \ \ \ \phi(b):=b^p+p\d b$$
is a ring homomorphism (equivalently a Frobenius lift); we say that $\d$ and $\phi$ are {\it attached} to each other. 
We view $p$-derivations as arithmetic analogues of derivations; cf. \cite{char, Jo}.
Then we view \ref{phipep} as analogues of Lie derivatives with respect to $p$-derivations. Similarly, if $X$ is a $p$-adic formal scheme over $A$, a $p$-{\it derivation} 
on $X$ (or an {\it arithmetic $\d^A$-flow} on $X$) is a map of sheaves of sets $\d=\d^X:\cO_X \ra \cO_X$ such that the map of sheaves of sets $\phi=\phi^X:\cO_X\ra \cO_X$, $\phi(b)=b^p+p\d b$, is a map of sheaves of rings (and hence induces a Frobenius lift  $\phi=\phi^X:X\ra X$). We again say that $\d$ and $\phi$ are {\it attached} to each other. 
As we will see later the above concept of {\it arithmetic $\d^A$-flow} is not flexible enough to accommodate some of the interesting examples of the theory; in the case of the Painlev\'{e} equations we will need a generalization  of the concept of {\it arithmetic $\d^A$-flow}
which will be referred to as {\it generalized arithmetic $\d^A$-flow}. 

 Let $\d$ be a $p$-derivation on some $p$-adically complete $p$-torsion free ring $B$. An element $c\in B$ is called a {\it $\d$-constant} 
if $\d c=0$. The set $B^{\d}\subset B$ of $\d$-constants is a multiplicative submonoid (but not a subring) of $B$. Let $\bZ[B^{\d}]$ be the subring of $B$ generated by $B^{\d}$ and let $\Sigma$ be the multiplicative system
$\Sigma:=B^{\times}\cap \bZ[B^{\d}]$.
An element  of $B$ is called 
a {\it pseudo-$\d$-constant} if it is a $p$-adic limit in $B$ of elements in the ring of fractions $\Sigma^{-1}\bZ[B^{\d}]$. So the 
set of pseudo-$\d$-constants in $B$ is a subring of $B$. 
One can easily check that 
if   $B/pB$ is perfect (i.e., the $p$-power Frobenius on $B/pB$ is surjective) then any element in $B$ is congruent mod $p$  to an element of $B^{\d}$ and, consequently, any element 
of $B$ is pseudo-$\d$-constant; in particular any element of $A$ is a pseudo-$\d$-constant.
Conversely, if  an element $b\in B$ is congruent mod $p$ to an element in $B^{\d}$ then $\d b$ is congruent mod $p$ to a $p$-th power in $B$.

One can introduce  arithmetic analogues of jet spaces as follows; cf. \cite{char}.
Let $X$ be a scheme of finite type over $A$ or the $p$-adic completion of such a scheme.
Say first that $X$ is affine,
$$X=Spec\ A[x]/(f)\ \ \text{or}\ \ \ X=Spf\ A[x]^{\widehat{\ }}/(f),$$
with $x$ and $f$ tuples.
 Then define the $p$-{\it jet spaces} of $X$ to be the $p$-adic formal schemes
\begin{equation}
\label{J}
J^n(X):=Spf\ A[x,x',...,x^{(n)}]^{\widehat{\ }}/(f,\d f, ... ,\d^nf)\end{equation}
where $x',...,x^{(n)}$ are new tuples of  indeterminates and $\d=\d^{\text{univ}}$ is the unique $p$-derivation on $A[x,x',...,x^{(n)},...]$ extending $\d^A$ and satisfying 
$\d x=x'$, ..., $\d x^{(n)}=x^{(n+1)}$,...
 We denote, as usual, by $\phi=\phi^{\text{univ}}$ the Frobenius lift attached to $\d^{\text{univ}}$; it induces ring homomorphisms 
 $$\phi=\phi^{\text{univ}}:\cO(J^n(X))\ra \cO(J^{n+1}(X)).$$
If $X$ is not necessarily affine then, again,  one defines $J^n(X)$ by gluing $J^n(X_i)$ where $X=\cup X_i$ is a Zariski affine open cover. (In the gluing process one uses the fact that we are dealing with formal schemes rather than schemes. There is a more global approach, avoiding gluing, that leads to functorially constructed algebraizations of our $p$-jet spaces; cf. \cite{borger1, borger2}. We will not need these algebraized $p$-jet spaces in what follows.) Then $(J^n(X))_{n\geq 0}$ has, again,  a structure of projective system of $p$-adic formal schemes depending functorially on $X$, with $J^0(X)=\widehat{X}$. 

If $G$ is a group in the category of schemes or $p$-adic formal schemes over $A$ then $(J^n(G))_{n\geq 0}$ has, again,  a structure of projective system of groups in the category of $p$-adic formal schemes; however, even if $G/A$ is smooth, the kernels 
of the homomorphisms $J^n(G)\ra J^{n-1}(G)$ are {\it generally not} isomorphic as groups to powers
of $\widehat{{\mathbb G}_a}$, although they are always isomorphic as formal schemes to some completed affine space $\widehat{{\mathbb A}^d}$. (By the way, these kernels are  commutative if and only if $G$ itself is commutative!) 

By an  {\it arithmetic differential equation} on $X$ we will understand a closed formal subscheme of some $J^n(X)$.
An {\it arithmetic  $\d^A$-flow} on $X$ will mean an arithmetic $\d^A$-flow on $\widehat{X}$.
To give an arithmetic $\d^A$-flow
on $\widehat{X}$  is equivalent to giving a section of the canonical projection $J^1(X)\ra \widehat{X}$, i.e. to giving a differential equation $Z\subset J^1(X)$ for which the projection $Z\ra \widehat{X}$ is an isomorphism. A {\it prime integral} for an arithmetic  $\d^A$-flow $\d^X$ is a 
$\d^X$-constant in $\cO(\widehat{X})$, i.e., a function $H\in \cO(\widehat{X})$ such that $\d^XH=0$.
For any $A$-point 
$P\in X(A)$, one defines the jets $J^n(P)\in J^n(X)(A)$ as the unique morphisms 
lifting $P$ that are  compatible with the actions of $\d^{\text{univ}}$ and $\d^A$. 
A {\it solution} (in $A$) for a differential equation $Z\subset J^n(X)$ is, again, an $A$-point 
$P\in X(A)$ such that $J^n(P)$ factors through $Z$. If $P$ is a solution to the differential equation defined by an arithmetic $\d^A$-flow and if $H$ is a prime integral for that arithmetic $\d^A$-flow
then, again, $\d^A(H(P))=0$; so, again, intuitively $H$ is ``constant along any solution".
 If $X$ is affine 
and $Z\subset J^1(X)$
is the arithmetic differential equation corresponding to an arithmetic  $\d^A$-flow $\d^X$ on $X$ then a point $P\in X(A)$ is a solution to $Z$ if and only if the ring homomorphism $P^*:\cO(\widehat{X})\ra A$ defined by $P$ 
satisfies $\d^A\circ P^*=P^*\circ \d^X$.

Let, in what follows, $S=Spec\ B$ be a smooth affine surface. Then, by the discussion in the previous subsection, we have a notion of {\it symplectic} form $\eta\in \Omega^2_{B/A}$ and associated {\it Poisson structure}, $\{\ ,\ \}_{\eta}$. For an arithmetic  $\d^A$-flow $\d=\d^{\widehat{B}}:\widehat{B}\ra \widehat{B}$ on $S$ the analogue of Lie derivatives will be the maps
\begin{equation}
\frac{\phi^*}{p^i}:\Omega^i_{\widehat{S}}\ra \Omega^i_{\widehat{S}},\ \ i=0,1,2.\end{equation}
At this point we would like to define what it means for an arithmetic  $\d^A$-flow on $S$ to be {\it Hamiltonian with respect to a symplectic form $\eta$ on $S$}.  One is tempted to make the following definition: an arithmetic $\d^A$-flow $\d=\d^S$ on $S$ is {\it Hamiltonian} with respect to the symplectic form $\eta$ on $S$
if
\begin{equation}
\label{Ham}
\frac{\phi^*}{p^2}\eta=\lambda\cdot \eta, \end{equation}
where $\lambda\in \cO(\widehat{S})$ is a pseudo-$\d$-constant. 
The concept we just defined is, however, 
not flexible enough to accommodate our examples. 
In particular, for the Painlev\'{e} equations one will 
need to replace arithmetic $\d^A$-flows with 
what we will call {\it generalized arithmetic $\d^A$-flows}; 
while for the Euler equations one will need  
to replace equality in \ref{Ham} by a congruence mod $p$. In view of the above
we will {\it not} adopt, in what follows, the above attempted definition of the {\it Hamiltonian} property but rather postpone the discussion of the Hamiltonian-related concepts to the next sections where the main examples of the theory are discussed; there we will encounter arithmetic analogues of the Hamiltonian property, canonical contact forms, canonical (generalized) arithmetic $\d^A$-flows, Euler-Lagrange forms, etc., each of which will be adapted to their specific context.

\section{Painlev\'{e} equations}

\subsection{The classical case} As a step in his proof of the Mordell conjecture over function fields of characteristic zero \cite{maninmordell} Manin introduced a differential algebraic map now referred to as the {\it Manin map}. A different, but equivalent, construction of this map was given in  \cite{annals}.
On the other hand Manin showed in \cite{maninmirror} (cf. also \cite{maninfrob}, p. 71)
how the Painlev\'{e} VI equation 
can be understood as a ``deformation of the Manin map"; he attributes this viewpoint to Fuchs.
In \cite{maninmirror} Manin also explained how this viewpoint leads to a Hamiltonian structure for the Painlev\'{e} VI equation. We quickly review here the ``deformation of the Manin map" interpretation of Painlev\'{e} VI in \cite{maninmirror}  and refer to \cite{maninmirror} for the Hamiltonian picture. 

Let $A$ be  the algebraic closure of a function field of one variable, 
$A=\overline{{\mathbb C}(t)}$, and $\d^A=d/dt$. 
Let ${\mathcal E}$ be an elliptic curve (i.e., a smooth projective curve of genus one) over $A$ which does not descend to ${\mathbb C}$.
Then the second jet space $J^2({\mathcal E})$ is easily seen to possess a non-zero group homomorphism of algebraic groups, unique up to multiplication by an element of $A$,
\begin{equation}
\label{mm}
\psi:J^2({\mathcal E})\ra {\mathbb G}_a,\end{equation}
into the additive group ${\mathbb G}_a$ over $A$. The map \ref{mm} is an incarnation
of the Manin map, as explained, in a more general setting, in \cite{annals}. Let us view
$\psi$ as an element of $\cO(J^2({\mathcal E}))$ and, for any open set $Y\subset {\mathcal E}$, let us view
$\cO(Y)$ and $\cO(J^2({\mathcal E}))$ as subrings of $\cO(J^2(Y))$ via pull-backs. 
Also let us recall that the classical Painlev\'{e} VI equation is a family, depending on $4$ parameters in ${\mathbb C}$, of differential equations.
Then Manin's analysis in \cite{maninmirror} shows  that each of the differential equations  in the  Painlev\'{e} VI family can be interpreted (in our language introduced above) as 
the closed subscheme $Z$ of $J^2(Y)$ defined by
\begin{equation}
\label{pain}
f:=\psi-r\in \cO(J^2(Y))
\end{equation}
where $Y$ is the complement in ${\mathcal E}$ of the set ${\mathcal E}[2]$ of $2$-torsion points in ${\mathcal E}(A)$
and $r\in \cO(Y)$ is an appropriate function. More precisely $r$ is a suitable  ${\mathbb C}$-linear combination of the $4$ translates, by the $4$  points in ${\mathcal E}[2]$, of the $y$-function on ${\mathcal E}\backslash {\mathcal E}[2]$ in a representation 
$${\mathcal E}\backslash {\mathcal E}[2]=Spec\ A[x,y,y^{-1}]/(y^2-(x^3+ax+b));$$  the 
complex coefficients of this linear combination are related to the $4$ complex parameters 
in the corresponding classical Painlev\'{e} equation. 
Moreover one can easily show that

\begin{thm}
For any function $r\in \cO(Y)$
the differential equation  $Z\subset J^2(Y)$ given by \ref{pain} defines a canonical $\d^A$-flow on $S:=J^1(Y)$.\end{thm}

In particular the equations in the Painlev\'{e} VI family are defined by canonical $\d^A$-flows on 
$J^1(Y)$. By the way notice that $J^1({\mathcal E})$, on which Painlev\'{e} equations ``live", is an
${\mathbb A}^1$-fibration over the elliptic curve ${\mathcal E}$; and  actually, over $Y$, this fibration is trivial, so we have an isomorphism 
$$J^1(Y)\simeq 
Y \times {\mathbb A}^1.$$
For details on the Hamiltonian picture we refer to \cite{maninmirror}. 

\subsection{The arithmetic case} The construction of the Manin map in \cite{annals}
was shown in \cite{char} to have an arithmetic analogue. Then, in \cite{BYM}, an arithmetic analogue of the Painlev\'{e} VI equation was introduced and a {\it Hamiltonian structure} was shown to exist for it. We explain this in what follows.

Let $A$ be a complete discrete valuation ring with maximal ideal generated by  $p$ and perfect residue field. Recall that elliptic curves over $A$ do not generally admit Frobenius lifts;
an elliptic curve that admits a Frobenius lift is automatically with complex multiplication.

\begin{thm} \label{charr} \cite{char}
 Let ${\mathcal E}$ be an elliptic curve over $A$ that  admits no  Frobenius lift.
There exists a non-zero group homomorphism, in the category of  $p$-adic formal schemes,
\begin{equation}
\label{amm}
\psi:J^2({\mathcal E})\ra \widehat{{\mathbb G}_a},\end{equation}
 which is unique up to multiplication by a constant in $A^{\times}$.\end{thm}

 We view \ref{amm} as an arithmetic analogue of the Manin map \ref{mm}.
Given an invertible $1$-form $\omega$ on ${\mathcal E}$ one can normalize $\psi$ with respect to $\omega$; we will need, and review, this normalization later.  The normalized $\psi$ can be referred to as the {\it  canonical $\d$-character} on ${\mathcal E}$; cf.  \cite{book}, Definition
7.24. One can view $\psi$ as an element of $\cO(J^2({\mathcal E}))$. 
By the way one has:

\begin{thm}\label{jeje}\cite{je}
Let ${\mathcal E}$ be an elliptic curve over $A$ that admits no Frobenius lift. Then the following hold:

1) $\cO(J^1({\mathcal E}))=R$.

2) $J^1({\mathcal E})$ admits no Frobenius lift.
\end{thm}

Assertion 1 in Theorem \ref{jeje} shows that ``order $2$ in Theorem \ref{charr} is optimal."
Assertion 2 in Theorem \ref{jeje} is equivalent to saying that the projection 
$$J^1(J^1({\mathcal E}))\ra J^1({\mathcal E})$$ 
does not admit a section in the category of $p$-adic formal schemes; equivalently, there is no arithmetic $\d^A$-flow on $J^1({\mathcal E})$! This justifies our generalization of the notion of arithmetic $\d^A$-flow  below.
To introduce this let us assume, for a moment that $Y$ is any smooth affine curve over $A$ and assume we are given 
 an arithmetic differential equation 
\begin{equation}
\label{mortify}
Z\subset \cO(J^r(Y)).
\end{equation}
  Then one can consider the module
\begin{equation}
\label{asthma77}
\Omega_Z = \lim_{\leftarrow} \Omega_{Z_n/A_n},
\end{equation}
$A_n:=A/p^nA$, $Z_n:=Z\otimes A_n$,
and the module
\begin{equation}
\label{asthma77}
\Omega_J = \lim_{\leftarrow} \Omega_{J_n/A_n},
\end{equation}
where $J:=J^r(Y)$.
On the other hand
put
\begin{equation}
\label{4.6}
\Omega^{\prime}_Z:= \frac{\Omega_{J}}{\langle I_Z \Omega_{J}, dI_Z\rangle}
\end{equation}
where $I_Z$ is the ideal of $Z$ in $J^r(Y)$.
Moreover define $\Omega_Z^{\prime i}$ to be the $i$-th wedge power $\wedge^i\Omega_Z'$.
Under  quite general hypotheses the modules \ref{maxim} and \ref{4.6} coincide; we will not discuss this here but, rather, refer to \cite{foundations}, Lemma 3.165.

 Going back to $Z$ as in \ref{mortify},
 for each $s\leq r$, there is a natural map
$$\pi_{r,s} :\, Z\to J^s(Y).$$
We also have natural maps 
$$
\frac{\phi^{\text{univ}*}}{p^i} :\, \Omega^i_{J^{r-1}(Y)}\to \Omega^i_{J^{r}(Y)},
$$
inducing maps which we will denote by
$$
\frac{\phi^*_Z}{p^i}:\Omega^i_{J^{1}(Y)}\to \Omega^{\prime i}_{Z}.
$$
We say that   $Z\subset J^2(Y)$   defines a {\it generalized 
$\delta$-flow} on $J^1(Y)$, if the
induced map
$$
\pi_{2,1}^*\Omega_{J^1(Y)}\to \Omega^{\prime}_{Z}
$$
is injective, and its cokernel is annihilated by a power of $p$.
Under quite general conditions, if $Z$ defines an arithmetic $\d^A$-flow on $J^1(Y)$ then $Z$  defines a generalized arithmetic $\d^A$-flow on $J^1(Y)$; again we will not need this so we will not discuss these conditions here; but see, again, \cite{foundations}, Lemma 3.165.

 Now let $S$ be a smooth surface   over $A$ or the $p$-adic completion of such a surface. 
 Recall that a {\it symplectic form} on $S$ is an invertible $2$-form on $X$ over $A$; a {\it contact form} on $S$ is a $1$-form on $X$ over $S$ such that $d\nu$ is symplectic; and for $S=J^1(Y)$ with
 $Y$  a smooth curve over $A$, a $1$-form $\nu$ on $S$ is called {\it canonical}
if $\nu=f\beta$, where $f\in \cO(S)$ and $\beta$ is an $1$-form lifted
from $Y$.

  Let  $Y$ be a smooth affine curve over $A$ and let $f\in \cO(J^2(Y))$ 
be a function whose zero locus defines  a generalized arithmetic $\d^A$-flow on $S:=J^1(Y)$. 
The respective generalized arithmetic $\d^A$-flow is called  {\it Hamiltonian} with respect to the symplectic form $\eta$ on $S$,
 if $$\frac{\phi^*_Z}{p}\eta =\lambda\cdot \eta$$
in $\Omega^{\prime 2}_{Z}$ for some $\lambda \in A$; note that any element in $A$, hence in particular $\lambda$, is a pseudo-$\d$-constant;
so the definition we just gave is a generalized version of the definition we proposed  in \ref{Ham}.
 Assume, moreover, that $\eta =d\nu$ for some canonical $1$-form  $\nu$ on $S$.
Then we call
\begin{equation}
\epsilon:= \frac{\phi^*_Z}{p}\nu-\lambda \nu\in \Omega^{\prime}_{Z}
\end{equation}
the {\it Euler-Lagrange form} attached to $\nu$.

Now let ${\mathcal E}$ be an elliptic curve over $A$ that does not admit a  Frobenius lift
and let $\psi\in \cO(J^2({\mathcal E}))$ be the canonical $\d$-character with respect to an invertible $1$-form $\omega$ on ${\mathcal E}$. Consider the symplectic form  
$$\eta=\omega\wedge \frac{\phi^{\text{univ}*}}{p}\omega$$ on $J^1({\mathcal E})$.
 Let $Y\subset {\mathcal E}$ be an affine open set possessing  an \'{e}tale coordinate; this latter condition is satisfied, for instance,  if $Y={\mathcal E}\backslash {\mathcal E}[2]$. By the way, notice that $J^1({\mathcal E})$ is an
$\widehat{{\mathbb A}^1}$-fibration over the elliptic curve ${\mathcal E}$; and actually this fibration is trivial over $Y$, hence we have an isomorphism of formal schemes,
$$J^1(Y)\simeq 
\widehat{Y} \times \widehat{{\mathbb A}^1}.$$
 
 \begin{thm} \label{bbb} \cite{BYM}
 
1) There exists a canonical  contact form $\nu$ on $S:=J^1(Y)$ such that $d\nu=\eta$.

2) For any function $r\in \cO(Y)$ the differential equation $Z \subset J^2(Y)$ given by the 
zero locus of the function  
$$f=\psi-\phi^{\text{univ}}(r)\in \cO(J^2(Y))$$  defines a generalized arithmetic $\d^A$-flow on $S$ which  is Hamiltonian with respect to $\eta$.
\end{thm} 

In particular the symplectic form $\eta$ is exact and  the Euler-Lagrange form  $\epsilon$ is closed.
 The function $f$ in assertion 2 of the Theorem is our analogue of
 the Painlev\'{e} VI equation. By the Theorem it 
 defines a generalized arithmetic $\delta$-flow on $S$; however, it does not define an arithmetic  $\delta$-flow on $S$  which is our motivation for generalizing the definition of arithmetic $\delta$-flow. 
 Note the discrepancy with the classical case coming from replacing $r$ by $\phi^{\text{univ}}(r)$ in the expression of $f$ in Theorem \ref{bbb}.
Another discrepancy comes from the absence, in the arithmetic setting,
of an analogue of the $4$ constant parameters in the classical Painlev\'{e} equations.

\section{Euler equations}

In Manin's picture \cite{maninmirror} we have just reviewed the Painlev\'{e} VI equation
``lives" on an ${\mathbb A}^1$-fibration over an elliptic curve. On the other hand, the {\it Euler equation}
describing the motion of a rigid body with a fixed point, which we are discussing next, ``lives" on an elliptic fibration over ${\mathbb A}^1$. This already suggests an analogy between the geometries underlying these differential equations and their arithmetic analogues; we will make such analogies/links more precise below.

\subsection{The classical case}
We begin by reviewing the classical Euler equations from a purely algebraic point of view. Let $A$ be
either a field or a discrete valuation ring and assume  $2$ is invertible in $A$.
Let $x_1,x_2,x_3$ and $z_1,z_2$ be variables and let $a_1,a_2,a_3\in A$ be such that 
$$(a_1-a_2)(a_2-a_3)(a_3-a_1)\in A^{\times}.$$
We consider the quadratic forms,
$$
H_1:=\sum_{i=1}^3 a_ix_i^2\in A[x_1,x_2,x_3],\ \ \ H_2:=\sum_{i=1}^3x_i^2\in A[x_1,x_2,x_3].$$
Also we consider the affine spaces 
$
{\mathbb A}^2=\operatorname{Spec}\ A[z_1,z_2]$, ${\mathbb A}^3=\operatorname{Spec}\ A[x_1,x_2,x_3]
$
and the morphism
${\mathcal H}:{\mathbb A}^3\ra {\mathbb A}^2$
defined by $z_1\mapsto H_1$, $z_2\mapsto H_2$. 
For $i=1,2,3$ denote by $Z_i\subset {\mathbb A}^3$ the $x_i$-{\it coordinate plane} and
let
$$
L_1=Z_2\cap Z_3,\ \ L_2=Z_3\cap Z_1,\ \ L_3=Z_1\cap Z_2
$$
 be the {\it $x_i$-coordinate axes}. 
 Then  ${\mathcal H}$ is smooth on the complement of 
 $L_1\cup L_2 \cup L_3$. 
 For any $A$-point $c=(c_1,c_2)\in A^2={\mathbb A}^2(A)$ we set
$$E_c:={\mathcal H}^{-1}(c)=\operatorname{Spec}\ A[x_1,x_2,x_3]/(H_1-c_1,H_2-c_2),$$
and we let $i_c:E_c\ra {\mathbb A}^3$ be the inclusion. 
Consider the polynomial
$$
 N(z_1,z_2)= \prod_{i=1}^3(z_1-a_iz_2)\in A[z_1,z_2].$$
Then, for 
 $N(c_1,c_2)\in A^{\times}$,
 $E_c$ is disjoint from $L_1\cup L_2 \cup L_3$ and, in particular,
 $E_c$ is smooth over $A$: it is an affine elliptic curve. 
 Moreover $E_c$ comes equipped with a global $1$-form
 given by
  \begin{equation}
  \label{the form on intersections of two quadrics}
  \omega_c=i_c^*\frac{dx_1}{(a_2-a_3)x_2x_3}=i_c^*\frac{dx_2}{(a_3-a_1)x_3x_1}=i_c^*\frac{dx_3}{(a_1-a_2)x_1x_2}.\end{equation}
If one considers  the  smooth projective model ${\mathcal E}_c$ of $E_c$ then
$\omega_c$ extends to an invertible $1$-form on the whole of ${\mathcal E}_c$.
 In the discussion below a certain plane quartic will play a role; let us review this next. Consider  two more indeterminates $x,y$,  and consider the  polynomial 
   \begin{equation}
   \label{quartic}
     F:=((a_2-a_3)x^2+z_1-a_2z_2)((a_3-a_1)x^2-z_1+a_1z_2)\in A[z_1,z_2][x] .
     \end{equation}
   For any $c=(c_1,c_2)\in A^2$ set
   $$
   E'_c:=\operatorname{Spec}\ A[x,y]/(y^2-F(c_1,c_2,x)).
   $$
   Then we have a morphism
    $\pi: E_c\ra E'_c$
     given by
     $ x\mapsto x_3$, $y\mapsto (a_1-a_2)x_1x_2$.
     If $N(c_1,c_2)\in A^{\times}$, $E'_c$ is smooth over $A$ and
     $$
      \pi^*(\frac{dx}{y})=i_c^*\frac{dx_3}{(a_1-a_2)x_1x_2}=\omega_c.$$
     For $A$ a perfect  field and   $c_1,c_2$ satisfying $N(c_1,c_2)\neq 0$  we have that $E'_c$ is a smooth plane curve. If ${\mathcal E}'_c$  is its smooth projective model then we have an induced  isogeny of elliptic curves, 
     ${\mathcal E}_c\ra {\mathcal E}'_c$.
     
     Assume now, until further notice, that $A$ is a field of characteristic zero (classically $A={\mathbb C}$, the complex field), viewed as equipped with the trivial derivation $\d^A=0$, and 
consider  the $A$-derivation $\d=\d^B$ on the polynomial ring $B=A[x_1,x_2,x_3]$  given by
\medskip
 \begin{equation}
 \label{Euler system}
 \d x_1  =  (a_2-a_3)x_2x_3,\ \ 
  \d x_2  =  (a_3-a_1)x_3x_1,\ \ 
   \d x_3  =  (a_1-a_2)x_1x_2.\end{equation}
   \medskip
   We  refer to the derivation $\d$ as   the {\it classical Euler flow} on ${\mathbb A}^3$. 
   
   For any $c=(c_1,c_2)\in A^2$ with $N(c_1,c_2)\neq 0$
   denote by $\d_c$ the derivation on 
   $\cO(E_c)$ induced by the derivation $\d$ on $B$. We have the following trivially checked classical fact:

   \begin{thm}
   \label{classical theorem}
  
  \ 
  
  1) $H_1$ and $H_2$ are prime integrals for the classical Euler flow, i.e., 
  $$\d H_1=\d H_2=0.$$

   2) For any $c=(c_1,c_2)\in A^2$ with $N(c_1,c_2)\neq 0$ 
   the  Lie derivative $\d_c$ on $\Omega^1_{\cO(E_c)/A}$ annihilates the  $1$-form  $\omega_c$ on $E_c$:
  $$
  \d_c \omega_c=0.
  $$\end{thm}
  
 Condition 2 can be viewed as a {\it linearization} condition for the $\d^A$-flow $\d_c$ on $E_c$.
 It is equivalent to $\d_c$ having an extension to a vector field on the compactification ${\mathcal E}_c$. It is the condition in 2 and {\it not} the ``extension to the compactification" property that will have an arithmetic analogue.
 
 The classical Euler flow fits into the Hamiltonian paradigm. We explain this in what follows. Since we will later need a discussion of these concepts in the arithmetic case as well
 we revert in what follows to the case  when $A$ is either a field or a discrete valuation ring.
 Let $c_2\in A^{\times}$ and set
$$S_{c_2}:=Spec\ A[x_1,x_2,x_3]/(H_2-c_2)\subset {\mathbb A}^3,$$
 the {\it sphere of radius $c_2^{1/2}$}. Then 
$S_{c_2}$ is the scheme theoretic pullback via ${\mathcal H}$ of the line in ${\mathbb A}^2$ defined by 
$z_2-c_2$ and hence
the map
$$H:S_{c_2}\ra {\mathbb A}^1=Spec\ A[z_1]$$
induced by $z_1\ra H_1$ is smooth above the complement of the closed subscheme defined by 
$$N(z_1,c_2)=(z_1-c_2a_1)(z_1-c_2a_2)(z_1-c_2a_3).$$ 
Now consider the $2$-forms
$$\ \ \eta_1=\frac{dx_2\wedge dx_3}{x_1},\ \ 
\eta_2=\frac{dx_3\wedge dx_1}{x_2},\ \ \eta_3=\frac{dx_1\wedge dx_2}{x_3}$$
defined on the complements in $S_{c_2}$ of $Z_1, Z_2, Z_3$, respectively.
These forms glue together defining a symplectic form $\eta_{c_2}$ on $S_{c_2}$. If one considers the Poisson structure $\{\ ,\ \}$ on $\cO({\mathbb A}^3)=A[x_1,x_2,x_3]$ defined by
$$\{x_1,x_2\}=x_3,\ \ \{x_2,x_3\}=x_1,\ \ \{x_3,x_1\}=x_2,$$
then $H_2$ is a {\it Casimir} i.e., $\{H_2,-\}=0$, so
this Poisson structure induces a Poisson structure $\{\ ,\ \}_{c_2}$ on each $\cO(S_{c_2})$.
On $\cO(S_{c_2})$ we have
$$\{x_1,x_2\}_{c_2}=\frac{dx_1\wedge dx_2}{\eta_{c_2}},\ \ \{x_2,x_3\}_{c_2}=\frac{dx_2\wedge dx_3}{\eta_{c_2}},\ \ \{x_3,x_1\}_{c_2}=\frac{dx_3\wedge dx_1}{\eta_{c_2}},$$
hence the Poisson structure $\{\ ,\ \}_{c_2}$ on $\cO(S_{c_2})$ 
coincides with the Poisson structure $\{\ ,\ \}_{\eta_{c_2}}$ on $\cO(S_{c_2})$ defined by the symplectic form
$\eta_{c_2}$ (because the two Poisson structures coincide on the generators $x_1,x_2,x_3$
of $\cO(S_{c_2})$). 
In other words $S_{c_2}$ are {\it symplectic leaves} for our Poisson structure on 
 $\cO({\mathbb A}^3)$,
with corresponding symplectic forms $\eta_{c_2}$. Furthermore, if $\d$ is the {\it classical Euler flow} \ref{Euler system}
then $\d$ induces a derivation $\d_{c_2}$ on each $\cO(S_{c_2})$ and the Lie derivative on $2$-forms,  
$$\d_{c_2}:\Omega^2_{\cO(S_{c_2})/A}\ra  \Omega^2_{\cO(S_{c_2})/A}$$
is trivially seen to satisfy
\begin{equation}
\label{lie}
\d_{c_2} \eta_{c_2}=0.\end{equation}
In other words we have:

\begin{thm}\label{other words}
The  $\d^A$-flow $\d_{c_2}$ on $S_{c_2}$ is symplectic with respect to  $\eta_{c_2}$.
\end{thm}

The link between the $2$-forms $\eta_{c_2}$ and the $1$-forms $\omega_c$ is as follows.
Consider the $1$-forms 
$$\omega_1=\frac{dx_1}{(a_2-a_3)x_2x_3},\ \ \omega_2=\frac{dx_2}{(a_3-a_1)x_3x_1},\ \ \omega_3=\frac{dx_3}{(a_1-a_2)x_1x_2}$$
defined on 
$$
S_{c_2}\backslash (Z_2\cup Z_3), \ \ \ S_{c_2}\backslash (Z_3\cup Z_1), \ \ \ S_{c_2}\backslash (Z_1\cup Z_2),$$
 respectively. Recall that 
for any $c=(c_1,c_2)$ with $N(c_1,c_2)\in A^{\times}$ the restrictions of $\omega_1,\omega_2,\omega_3$ to $E_c$ glue to give the form $\omega_c$ on $E_c$.
A trivial computation then gives the following equalities of $2$-forms on $S_{c_2}\backslash (Z_1\cup Z_2\cup Z_3)$ which will play a role later:
\begin{equation}
\label{fiona}
\eta_{c_2}=-d H_1 \wedge \omega_1=-d H_1 \wedge \omega_2=-d H_1 \wedge \omega_3.\ \ 
\end{equation}
By the way, the equalities \ref{fiona} imply that for all $c=(c_1,c_2)$ with $N(c_1,c_2)\in A^{\times}$ 
the form $\omega_c$ on $E_c$ satisfies
\begin{equation}
\label{PR}
\omega_c=-P.R.\left(\frac{\eta_{c_2}}{H_1-c_1}\right),
\end{equation}
where 
$$P.R.:\Omega^2_{S_{c_2}/A}(E_c)\ra \Omega^1_{E_c/A}$$
is the Poincar\'{e} residue map \cite{GH}, p. 147; we will not need this interpretation in what follows.

\subsection{The arithmetic case}
In what follows $A$ is a complete discrete valuation ring with maximal ideal generated by an odd prime $p$ and perfect residue field $k=A/pA$.
Let $F\in A[z_1,z_2][x]$ be the polynomial in \ref{quartic}.
Define the {\it Hasse invariant}  to be the coefficient $A_{p-1}\in A[z_1,z_2]$ of $x^{p-1}$ in the polynomial
$F^{\frac{p-1}{2}}.$ In addition to the quantities defined in the previous
subsection we also consider the following polynomial
$$
  Q:=x_1x_2\cdot N(H_1,H_2)\cdot A_{p-1}(H_1,H_2)\in A[x_1,x_2,x_3],
 $$
 and the open subscheme of ${\mathbb A}^3$ defined by
 $$
  X=Spec\ A[x_1,x_2,x_3][1/Q].$$
     Assume in addition that $c=(c_1,c_2)\in A^2$ satisfies
  $$
  \d c_1=\d c_2=0\ \ \ \text{and}\ \ \ N(c_1,c_2)\cdot A_{p-1}(c_1,c_2)\in A^{\times}
  $$
    and let $\d^{X}$ be any arithmetic $\d^A$-flow on $\widehat{X}$ satisfying 
    $\d^X H_1=\d^X H_2=0$. Then the Frobenius lift 
  $\phi^{X}$ on $\cO(\widehat{X})$ induces a Frobenius lift
  $\phi_c:=\phi^{E^0_c}$  on 
  $\widehat{E^0_c}$
  where $E^0_c$ is the open set of $E_c$ 
  given by
  $E^0_c  :=  E_c\cap X$.
 We refer to $\phi_c$ as the {\it Frobenius lift} on $\widehat{E^0_c}$ attached to $\d^X$.
    On the other hand,  the  global $1$-form $\omega_c$ in
    \ref{the form on intersections of two quadrics}  restricted to $E^0_c$ will be referred to as the {\it canonical} $1$-form on $E^0_c$ and will still be denoted  by $\omega_c$.
    
   \medskip
   
    The following provides  an arithmetic analogue of the classical Euler flow: 
  assertions 1 and 2   below are arithmetic analogues of assertions 1 and 2 in Theorem \ref{classical theorem} respectively.
  
  \begin{thm}
  \label{linearization theorem}
 \cite{euler}
 There exists  an arithmetic $\d^A$-flow $\d^X$ on $\widehat{X}$ such that:
 
 1) $H_1$ and $H_2$ are prime integrals for $\d^X$, i.e., the following holds in $\cO(\widehat{X})$:
  $$
   \d^XH_1=\d^XH_2=0;$$
   
   2) For
 any point $c=(c_1,c_2)\in A^2$ with
 $$
 \d c_1=\d c_2=0,\ \ \ \text{and}\ \ \ N(c_1,c_2)\cdot A_{p-1}(c_1,c_2)\in A^{\times}$$
  the  Frobenius lift $\phi_c$ on $\widehat{E_c^0}$ attached to $\d^X$ and the canonical $1$-form $\omega_c$ on $E_c^0$ satisfy the following congruence in $\Omega^1_{\widehat{E^0_c}}$: 
  $$
  \frac{\phi_c^*}{p}\omega_c\equiv A_{p-1}(c_1,c_2)^{-1} \cdot \omega_c\ \ \ \text{mod}\ \ \ p.$$
\end{thm}

By the way, 
one can ask if the open set $X$ in Theorem \ref{linearization theorem} can be taken to be the whole of ${\mathbb A}^3$. 
In contrast with the classical case (Theorem \ref{classical theorem}), 
the answer to this is negative; indeed we have the following ``singularity theorem":

\begin{thm}
\label{singularity}
\cite{euler}
If $X\subset {\mathbb A}^3$ is an open set such that $\widehat{X}$ possesses an arithmetic $\d^A$-flow
$\d^X$ with $
   \d^XH_1=\d^XH_2=0$ and if $\d a_i\in A^{\times}$ for some $i\in\{1,2,3,\}$ then $\widehat{X}$ cannot meet the coordinate axis $\widehat{L_i}$.
\end{thm}

Another question one can ask is whether it is possible to extend the Frobenius lifts
$\phi_c$ in Theorem \ref{linearization theorem}
to the compactifications ${\mathcal E}_c$ of $E_c$. In contrast with the classical case (Theorem \ref{classical theorem}), 
the answer to this is, again,  negative, cf. Theorem \ref{masa} below.
For this theorem we fix, for every rational prime $p$, a complete discrete valuation ring $R_p$ with maximal ideal generated by $p$ and algebraically closed residue field. We also fix a number field $F$, with ring of integers $\cO_F$, a rational  integer $M$, and, for each $p>>0$ we fix an embedding of $\cO_F[1/M]$  into $R_p$.

 \begin{thm}
 \label{masa}
    \cite{canonical}
    Let $a_1,a_2,a_3\in \cO_F[1/M]$. Then, if $p>>0$, there is no triple $(K,X,\phi^X)$ with

  $\bullet$ $K\in \cO(\widehat{{\mathbb A}^2})=R_p[z_1,z_2]^{\widehat{\ }}$,  $K\not\equiv 0$ mod $p$,
    
  $\bullet$  
  $X\subset {\mathbb A}^3$ an open set over $R_p$, 
  
  $\bullet$ $\phi^X$   a Frobenius lift  on  $\widehat{X}$,

    \noindent   satisfying
    the following two conditions: 
    
    1) $H_1$ and $H_2$ are prime integrals for the arithmetic $\d^A$-flow $\d^X$ attached to $\phi^X$, i.e., the following holds in $\cO(\widehat{X})$:
  $$
   \d^XH_1=\d^XH_2=0;$$
     
 2) for all $c\in R_p^2$ with $\d c=0$, $N(c)K(c)\in R_p^{\times}$, one has
 $\widehat{E}_c\cap \widehat{X}\neq \emptyset$ and
 $$\textit{$\phi_c$ extends to an endomorphism of 
    the compactification ${\mathcal E}_c$ of $E_c$,}
   $$
  where   $\phi_c$ is  the Frobenius lift on  $\widehat{E}_c\cap \widehat{X}$ induced by $\phi^X$.
    \end{thm}

Interestingly, the proof of Theorem \ref{masa} is based on a variant of a Diophantine  result
in \cite{local} which, in its turn, is proved using, again,  arguments involving arithmetic differential equations.

\medskip

In what follows we use Theorem \ref{linearization theorem} to derive an arithmetic analogue of the Hamiltonian picture; cf. Theorem \ref{other words}.

Let $c_2\in A^{\times}$ be such that $\d c_2=0$. Then the Frobenius lift $\phi^X$
attached to the arithmetic $\d^A$-flow $\d^X$ in Theorem \ref{linearization theorem} induces a Frobenius lift $\phi_{c_2}$ on $\widehat{S^0_{c_2}}$ where $S^0_{c_2}:=S_{c_2}\cap X$.
 Recall that it follows from equations 6.1 and 6.2 in \cite{euler} that the function $A_{p-1}(H_1,c_2)$ is invertible on $\widehat{S^0_{c_2}}$. 
 Set
 $$\lambda:=\frac{H_1^{p-1}}{A_{p-1}(H_1,c_2)}\in \cO(\widehat{S_{c_2}^0})$$
 and note that $\lambda$  is a pseudo-$\d$-constant in
$\cO(\widehat{S_{c_2}^0})$ because $H_1$ is a  $\d$-constant and all elements of $A$ are pseudo-$\d$-constants. 
 We will prove the following result which can be interpreted as a relaxation of the condition defining the {\it Hamiltonian} property in \ref{Ham}:

\begin{theorem}
\label{new1} The following holds in $\Omega^2_{\widehat{S_{c_2}^0}}$:
$$\frac{\phi_{c_2}^*}{p^2}\eta_{c_2}\equiv \lambda \eta_{c_2}\ \ \ 
\text{mod}\ \ \ p.$$
\end{theorem}

 {\it Proof}. 
 Consider the form $\theta$ on $\widehat{S^0_{c_2}}$ defined by
$$ \theta:=\frac{\phi_{c_2}^*}{p^2}\eta_{c_2}-\frac{H_1^{p-1}}{A_{p-1}(H_1,c_2)}\eta_{c_2}.
$$
By \ref{fiona} and $\phi_{c_2}(H_1)=H_1^p$ we have
$$
\theta = -H_1^{p-1}d H_1\wedge  \frac{\phi_{c_2}^*}{p}\omega_1 +\frac{H_1^{p-1}}{A_{p-1}(H_1,c_2)}dH_1\wedge \omega_1= - H_1^{p-1}dH_1 \wedge \beta
$$
where 
$$\beta:=\frac{\phi_{c_2}^*}{p}\omega_1- \frac{1}{A_{p-1}(H_1,c_2)} \omega_1
$$
Let $i_c:E^0_c=E_c\cap X\ra S^0_{c_2}$ be the inclusion. Then, if $\d c_1=0$,
$N(c_1,c_2)\in A^{\times}$, $A_{p-1}(c_1,c_2)\in A^{\times}$, by Theorem \ref{linearization theorem},
$$i_c^*\beta= \frac{\phi_{c}^*}{p}\omega_c- \frac{1}{A_{p-1}(c_1,c_2)} \omega_c
\equiv 0\ \ \ \text{mod}\ \ \ p.$$
Let us denote by an upper bar the operation of reduction mod $p$. Since any element in $k$ can be lifted to an element $c_1$ of $A$ killed by $\d$ (this lift is the Teichm\"{u}ller lift) it follows that 
$$\overline{i}^*_{\overline{c}} \overline{\beta}=0$$
for all except finitely many $\overline{c}_1\in k$, where $\overline{i}_{\overline{c}}:\overline{E^0_c}\ra \overline{S^0_{c_2}}$ is the inclusion. 
Recall  from \cite{euler} that $H_1,H_2,x_3$ are \'{e}tale coordinates on $X$; so $H_1,x_3$ are \'{e}tale coordinates on $S_{c_2}^0$.
Write
$$\beta=b_1dH_1+b_2dx_3,$$
 on  $\widehat{S^0_{c_2}}$, with $b_1,b_2\in \cO(\widehat{S^0_{c_2}})$.
Since
$$\overline{i}^*_{\overline{c}} \overline{\beta}=\overline{i}^*_{\overline{c}}\overline{b}_1 \cdot d \overline{c_1}+\overline{i}^*_{\overline{c}}\overline{b}_2 \cdot d\overline{x_3}=\overline{i}^*_{\overline{c}}\overline{b}_2 \cdot d\overline{x_3}$$
it follows that 
$$\overline{i}^*_{\overline{c}}\overline{b}_2 =0.$$
Since this is true for all except finitely many $\overline{c}_1$ it follows that $\overline{b_2}=0$.
But then
$$\overline{\theta}=- \overline{H_1}^{p-1}d\overline{H_1} \wedge \overline{\beta}=
- \overline{H_1}^{p-1}d\overline{H_1} \wedge \overline{b}_2 d\overline{x_3}=0.$$
We conclude that
the congruence in the statement of the Theorem holds on an open set of $\widehat{S_{c_2}^0}$ and hence on the whole of $\widehat{S_{c_2}^0}$.\qed

\medskip

We next deduce a result that establishes a {\it link} between the Painlev\'{e} paradigm and the Euler paradigm. Assume 
we are under the hypotheses and notation of Theorem \ref{linearization theorem}, with
$a_1,a_2,a_3\in \bZ_p$. (So morally we are in the ``Euler paradigm".)
 Assume moreover that $c_1,c_2$ in assertion 2 of that Theorem belong to $\bZ_p$ and are such that
${\mathcal E}_c$ does not have a Frobenius lift. Examples of this situation are abundant; cf. the last Remark in \cite{canonical}. 
Let, furthermore,  $\phi_c:\widehat{E^0_c}\ra \widehat{E^0_c}$ be as in Theorem \ref{linearization theorem} and let $\sigma^n_c:\widehat{E^0_c}\ra J^n(E^0_c)$ be the sections of the projections
$J^n(E^0_c)\ra \widehat{E^0_c}$
induced by $\phi_c$.
On the other hand let  $\psi_c\in J^2({\mathcal E}_c)$ be the canonical $\d$-character. (The latter belongs, as we saw, to the ``Painlev\'{e} paradigm".) 
 Assume, for simplicity, that the field  $k:=A/pA$ is algebraically closed
 and let $K_c$ be function field of ${\mathcal E}_c\otimes k$.
  We will prove:

\begin{thm}
\label{new2}
The image of $\sigma_c^{2*} \psi_c$ in $K_c$ is a $p$-th power in $K_c$.\end{thm}

{\it Proof}.  Recall by \cite{book}, Corollary 7.28, that
\begin{equation}
\label{opop}
d\psi_c=\lambda_2  \left(\frac{\phi_c^{\text{univ}*}}{p}\right)^2\omega_c+\lambda_1\frac{\phi_c^{\text{univ}*}}{p}\omega_c+\lambda_0\omega_c\end{equation}
in $\Omega_{J^2({\mathcal E}_c)}$ where $\phi_c^{\text{univ}}:J^n({\mathcal E}_c)\ra J^{n-1}({\mathcal E}_c)$ are the  universal Frobenius lifts and $\lambda_i\in A$, $\lambda_2=p$. By the way, the above holds without the assumption that $a_i, c_j\in \bZ_p$; also the equality $\lambda_2=p$ is precisely the definition of $\psi$ being {\it normalized} with respect to $\omega_c$. With the additional assumption that $a_i, c_j\in \bZ_p$ we have that ${\mathcal E}_c$ descends to  an elliptic curve
${\mathcal E}_{c/\bZ_p}$ over $\bZ_p$; then,
by \cite{frob}, Theorem 1.10,  and \cite{book}, Theorem 7.22, we actually also have
$$\lambda_1=-a_p,\ \ \ \lambda_0=1,$$
where $a_p$ in the trace of Frobenius acting on the reduction mod $p$ of  ${\mathcal E}_{c/\bZ_p}$. 
Similarly ${\mathcal E}'_c$ descends to an elliptic curve ${\mathcal E}'_{c/\bZ_p}$ over $\bZ_p$.
 Since there is a separable  isogeny between ${\mathcal E}'_{c/\bZ_p}$ 
and ${\mathcal E}_{c/\bZ_p}$  it follows that $a_p$ is also the trace of Frobenius acting on the reduction mod $p$ of ${\mathcal E}'_{c/\bZ_p}$.
On the other hand, by \cite{silverman}, p. 141-142,
$$a_p\equiv 1-|{\mathcal E}'_{c/\bZ_p}({\mathbb F}_p)|\ \ \ \text{mod}\ \ \ p.$$
Now by \cite{euler}, Lemma 5.2, 
$$|{\mathcal E}'_{c/\bZ_p}({\mathbb F}_p)|\equiv 1-A_{p-1}(c_1,c_2)\ \ \ \text{mod}\ \ \ p.$$
It follows that 
$$\lambda_1\equiv -A_{p-1}(c_1,c_2)\ \ \ \text{mod}\ \ \ p.$$
Let us view the maps $\cO(J^n({\mathcal E}_c))\ra \cO(J^{n+1}({\mathcal E}_c))$ induced by the natural projections 
as inclusions. Then $\phi_c$ equals the composition 
$\phi_c^{\text{univ}}\circ \sigma_c^1$; hence 
$$\sigma_c^{2*} \phi^{\text{univ}*}=\sigma_c^{1*} \phi^{\text{univ}*}=\phi_c^*.$$
 (Note, by the way,  that $\phi_c^2$ is not equal to $(\phi_c^{\text{univ}})^2\circ \sigma_c^2$!)
Taking $\sigma_c^{2*}$ in \ref{opop} we get
$$d (\sigma_c^{2*}\psi_c)=\sigma_c^{2*} d\psi_c\equiv -A_{p-1}(c_1,c_2)\frac{\phi_c^*}{p}\omega_c+\omega_c\equiv 0\ \ \ \text{mod}\ \ \ p.$$
If $K_c$ is the function field of ${\mathcal E}_c\otimes k$
we have $K_c=k(x,\gamma)$ with $x$ a variable and $\gamma$ quadratic over $k(x)$. Since $k(x,\gamma)=k(x,\gamma^p)$ we may write 
$$\overline{\sigma_c^{2*}\psi_c}=u+v\gamma^p\in K_c,\ \ \ u,v\in k(x)$$
hence
$$0=d (\overline{\sigma_c^{2*}\psi_c})=\left(\frac{du}{dx}+\frac{dv}{dx} \gamma^p\right)dx\in \Omega_{K_c/k}=K_cdx,$$
hence $\frac{du}{dx}=\frac{dv}{dx}=0$ which implies that $u,v\in k(x^p)$ as one can see by 
considering the simple fraction decomposition of $u$ and $v$. Consequently
$\overline{\sigma_c^{2*}\psi_c}\in K_c^p$. \qed

\section {Lax equations}

\subsection{The classical case}
Let $A$ be a Noetherian ring, let $B:=A[x_1,...,x_N]$ be a polynomial ring, and consider the affine space
${\mathbb A}^n=Spec\ B$.
Let $L$ be an
 $A$-Lie algebra, free as an $A$-module, with basis $e_1,...,e_N$, and write
$$[e_i,e_j]=\sum_k c_{ijk}e_k,\ \ \ c_{ijk}\in A.$$
Then there is a unique Poisson structure $\{\ ,\ \}$ on $B$, 
(or on ${\mathbb A}^N=Spec\ B$)
called the {\it Lie-Poisson} structure attached to $(L,(e_i))$, such that
$$[x_i,x_j]=\sum_k c_{ijk}x_k.$$

In particular we may consider the variables $x_1,...,x_N$, with $N=n^2$,  to be the entries of a matrix of indeterminates $x=(x_{ij})$, we may consider
 the affine space
${\mathfrak g}:={\mathbb A}^{n^2}=Spec\ B$, $B:=A[x]$, and we may 
consider the Lie algebra  $L:={\mathfrak g}(A)$
of $n\times n$ matrices with coefficients in $A$, with respect to the commutator,
with basis $e_{ij}$ the matrices that have $1$ on position $(i,j)$ and $0$ everywhere else. One may consider then 
the Lie-Poisson structure $\{\ ,\ \}$ on $B$ (equivalently on ${\mathfrak g}$)
attached to $(L,(e_{ij}))$.

On the other hand
 let $\d=\d^A$ be a derivation on $A$,  and let  $\d=\d^{\mathfrak g}$ be a $\d^A$-flow on ${\mathfrak g}$, i.e.  $\d^{\mathfrak g}$ is a derivation on $B=A[x]$, 
extending $\d^A$.
Say that  $\d=\d^{\mathfrak g}$  is a {\it Lax $\d^A$-flow} if we have an equality of matrices with $B$-coefficients
$$\d x=[M,x]$$
 for some matrix $M=(m_{ij})$ with $B$-coefficients, i.e., 
 $$\d x_{ij}=\sum_k (m_{ik}x_{kj}-x_{ik}m_{kj}).$$

It is trivial to check that any $\d^A$-flow on $A[x]$ that is Hamiltonian with respect to the Lie-Poisson structure on $A[x]$ is a Lax $\d^A$-flow on ${\mathfrak g}$: if the Hamiltonian is
$H$ then $M$ can be taken to be the matrix $\frac{\partial H}{\partial x}:=\left(\frac{\partial H}{\partial x_{ij}}\right)$. 
Let us say that a Lax $\d^A$-flow on ${\mathfrak g}$ is {\it Hamiltonian} 
(or more accurately {\it Poisson-Hamiltonian}) if it  is Hamiltonian with respect to the Lie-Poisson structure on $A[x]$, equivalently if
$$\d x=[\frac{\partial H}{\partial x},x]$$
for some $H\in B$.

On the other hand, assuming for simplicity that $A$ is an algebraically closed field of characteristic zero,  any Lax $\d^A$-flow $\d$ is {\it isospectral}, by which we understand that:

\begin{thm}\label{par}
 The following diagram is commutative:
$$
\begin{array}{rcl}
B & \stackrel{\d}{\longrightarrow} & B\\
{\mathcal P} \uparrow & \  & \uparrow {\mathcal P}\\
A[z] &\stackrel{\d_0}{\longrightarrow} & A[z]\end{array}
 $$
 where $A[z]=A[z_1,...,z_n]$, a polynomial ring in $n$ variables, 
  $\d_0:A[z]\ra A[z]$ is the unique derivations extending $\d^A$ with
  $\d_0z_j=0$, and 
  ${\mathcal P}:A[z]\ra B$ is the $A$-algebra homomorphism with ${\mathcal P}(z_j)={\mathcal P}_j(x)$, where
  $$\det(s\cdot 1_n-x)=\sum_{j=0}^n (-1)^j{\mathcal P}_j(x)s^{n-j}.$$ \end{thm}

 In the above $1_n$ is the identity matrix, $s$ is a variable, 
 ${\mathcal P}_0=1$, and, for $j=1,...,n$,  ${\mathcal P}_j(x)$ are, of course, the coefficients of the characteristic polynomial of $x$:
 $${\mathcal P}_1(x)=\text{tr}(x),...,{\mathcal P}_n(x)=\det(x).$$
The commutativity of the above diagram implies $\d^A({\mathcal P}_j(x))=0$, i.e., ${\mathcal P}_j(x)$ are prime integrals
  for any Lax $\d^A$-flow. This implies that the characteristic polynomial of any solution to a Lax $\d^A$-flow  has $\d$-constant coefficients;  equivalently, the spectrum of any solution consists of $\d$-constants. (This equivalence will fail in the arithmetic case.) 

\subsection{The arithmetic case}
As usual we consider a complete discrete valuation ring $A$ with maximal ideal generated by an odd  rational  prime $p$ and perfect residue field and we view $A$ as equipped with its unique $p$-derivation $\d=\d^A$.
There are two arithmetic analogues of Lax $\d^A$-flows: one for which the characteristic polynomial of any solution
has $\d$-constant coefficients; and another one for which the solutions have $\d$-constant spectrum. These two conditions are not equivalent because, for a monic polynomial 
$$\sum_{j=0}^n a_j s^j=\prod_{j=1}^n (x-r_j) \in A[s],$$
 with all its roots $r_j$ in $A$, the condition that $\d a_j=0$ for all $j$ is {\it not} equivalent to the condition $\d r_j=0$ for all $j$;
these two conditions are equivalent for $\d$ a derivation on a field of characteristic zero but {\it not} for $\d$ our $p$-derivation on $A$. In what follows we explain these two analogues of Lax equations following \cite{foundations}. 

First let $T\subset G:=GL_n$ be the diagonal maximal torus,
$$T=Spec\ A[t_1,t_1,^{-1},...,t_n,t_n^{-1}],\ \ \ G=Spec\ A[x,\det(x)^{-1}],$$
with embedding given by $x_{jj}\mapsto t_j$ and $x_{ij}\mapsto 0$ for $j\neq i$, 
 and consider the map
$${\mathcal C}:T\times G\ra G,\ \ {\mathcal C}(h,g)=g^{-1}hg.$$
We have the following:

\begin{thm}\label{t*} \cite{foundations}
There exists an open set $G^*$ of $G=GL_n$  and  a unique Frobenius lift $\phi^{G^*}$ on $\widehat{G^*}$ such that the following diagram is commutative:
 $$ \begin{array}{rcl}
  \widehat{T^*}\times \widehat{G} & \stackrel{\phi^{T^*}_0\times \phi^G_{0}}{\longrightarrow} & \widehat{T^*}\times \widehat{G}\\
  {\mathcal C}\downarrow & \  & \downarrow {\mathcal C}\\
  \widehat{G^*} & \stackrel{\phi^{G^*}}{\longrightarrow} & \widehat{G^*}\end{array}
  $$
  where $T^*:=T\cap G^*$, ${\mathcal C}(T^*\times G)\subset G^*$, $\phi^{T^*}_0$ is induced by the unique Frobenius lift on $T$ that sends $t_j\mapsto t_j^p$, and $\phi^G_{0}$ is the Frobenius lift on $\widehat{G}$ that sends $x_{ij}\mapsto x_{ij}^p$.
 \end{thm}

Cf. \cite{foundations}, Theorem 4.50.
 By the way, 
in contrast with the classical case (Theorem \ref{par}), and in analogy with the arithmetic Euler paradigm (Theorem \ref{singularity}) we have the following ``singularity theorem":

\begin{thm} \cite{foundations}
For $n\geq 3$,  $G^*$ in Theorem \ref{t*} cannot be taken to be the whole of $G$.
\end{thm}

Cf. \cite{foundations}, Theorem 4.54.  
Also, it was shown in \cite{foundations}, Theorem 4.60, that any solution to the arithmetic $\d^A$-flow $\d^{G^*}$ attached to $\phi^{G^*}$, with spectrum contained in $A$, 
has the property that its spectrum consists of $\d^A$-constants. 
More generally, the same property holds if one replaces the Frobenius lift 
$\phi^{G^*}$ on $\widehat{G^*}$
by any Frobenius lift $\phi^{G^*(\alpha)}$ on $\widehat{G^*}$ that is {\it conjugate} to 
$\phi^{G^*}$ in the sense that
$$
\phi^{G^*(\alpha)}(x) :=\epsilon(x)^{-1} \cdot \phi^{G^*}(x) \cdot \epsilon(x),$$
where $\epsilon(x)=1+p\alpha(x)$, $\alpha(x)$ any $n\times n$ matrix with coefficients in $\cO(\widehat{G^*})$. By the way, for $\alpha$ with coefficients in $A$ (rather than in
$\cO(\widehat{G^*})$) the arithmetic $\d^A$-flows corresponding to 
$\phi^{G^*}$ and $\phi^{G^*(\alpha)}$ are {\it linear} with respect to each other in the sense of \cite{foundations}; we will not review this concept of linearity here.
We refer to loc. cit.  for details.

On the other hand we have:

\begin{thm}\label{t**} \cite{foundations}
There exists an open set $G^{**}$ of $G=GL_n$  and  a  Frobenius lift $\phi^{G^{**}}$ on $\widehat{G^{**}}$ such that the following diagram is commutative:
 $$
\begin{array}{rcl}
\widehat{G^{**}}& \stackrel{\phi^{G^{**}}}{\longrightarrow} & \widehat{G^{**}}\\
{\mathcal P} \downarrow & \  & \downarrow {\mathcal P}\\
\widehat{{\mathbb A}^n} &\stackrel{\phi^{{\mathbb A}^n}_0}{\longrightarrow} & \widehat{{\mathbb A}^n},\end{array}
 $$
 where $\phi^{{\mathbb A}^n}_0$ is induced by the unique Frobenius lift on ${\mathbb A}^n=Spec\ A[z_1,...,z_n]$ which sends $z_j\mapsto z_j^p$.
\end{thm}

Cf. \cite{foundations}, Theorem 4.56.
The polynomials ${\mathcal P}_j(x)$
are then prime integrals for  the arithmetic $\d^A$-flow $\d^{G^{**}}$ attached to $\phi^{G^{**}}$. In particular
the characteristic polynomial of any solution to the arithmetic $\d^A$-flow $\d^{G^{**}}$  are $\d^A$-constant. More generally, the same property holds if one replaces the Frobenius lift 
$\phi^{G^{**}}$ on $\widehat{G^{**}}$
by any Frobenius lift $\phi^{G^{**}(\alpha)}$ on $\widehat{G^{**}}$ that is {\it conjugate} to 
$\phi^{G^{**}}$ in the same sense as before, namely that 
$$
\phi^{G^{**}(\alpha)}(x) :=\epsilon(x)^{-1} \cdot \phi^{G^{**}}(x) \cdot \epsilon(x),$$
where $\epsilon(x)=1+p\alpha(x)$, $\alpha(x)$ any $n\times n$ matrix with coefficients in $\cO(\widehat{G^{**}})$. Again, for $\alpha$ with coefficients in $A$ (rather than in
$\cO(\widehat{G^{**}})$) the arithmetic $\d^A$-flows corresponding to 
$\phi^{G^{**}}$ and $\phi^{G^{**}(\alpha)}$ are {\it linear} with respect to each other in the sense of \cite{foundations}. We refer to loc. cit. for details.
The $\phi^{G^{**}}$ in  Theorem \ref{t**} is not unique; one can further subject it
to appropriate constraints that make it unique; we will not go into this here. 

In view of the above mentioned ``isospectrality-type" properties for the arithmetic $\d^A$-flows $\d^{G^*}$ and $\d^{G^{**}}$ in Theorems \ref{t*} and \ref{t**} respectively one may see these arithmetic flows as analogues of the classical Lax $\d^A$-flows. One is then tempted to ask for
an arithmetic analogue of the condition that a Lax $\d^A$-flow be   Poisson-Hamiltonian, i.e., an arithmetic analogue of the condition that the matrix $M$ in the classical equation $\d x=[M,x]$ is of the form $M=\frac{\partial H}{\partial x}$ for some  $H\in B$. The matrix $M$ itself does not have an obvious arithmetic analogue so the problem needs to be approached on a more conceptual level.

\end{document}